\documentclass[11pt,a4paper]{article}
\usepackage[margin=2cm]{geometry}
\usepackage[pdftex]{graphicx}
\DeclareGraphicsExtensions{.pdf,.png}

\usepackage{amsmath}
\usepackage{amsfonts}
\usepackage{mathrsfs}
\usepackage[utf8]{inputenc}

\newcommand{\const}{\operatorname{const.}}

\usepackage[T1]{fontenc}
\usepackage{amsthm}
\usepackage{amssymb}
\usepackage{cite}

\newcommand\gap{\operatorname{gap}}

\catcode`_=12
\begingroup\lccode`~=`_\lowercase{\endgroup\let~\sb}
\mathcode`_="8000

\title{
Gaps in the fractional parts of square roots
}
\author{Simon Čopar\\
  Faculty of Mathematics and Physics, University of Ljubljana\\
  Jadranska 19, Slovenia\\
  \texttt{simon.copar@fmf.uni-lj.si}
}
\date{\today}

\begin{document}
\maketitle

\begin{abstract}
  Fractional parts of the first $N$ natural numbers fill
  the unit interval with asymptotically uniform density.
  However, the gaps around rational points shrink at an
  asymptotically lower rate $N^{-1/2}$,
  and their widths scale with the Thomae (``popcorn'') function.
  This curious connection is derived and related geometrically
  to shadow pattern in the Euclid's orchard.
  Generalized cases of higher radicals and their convergence rates, are also investigated.
\end{abstract}

\section{Introduction}

Fractional parts of real number sequences incrementally fill the unit
interval with values in a way that strongly depends on the studied
sequence. The distributions of values and the gaps between them in the
limit of an infinite sequence have been studied from different
perspectives. A lot of effort has been put into researching the
distribution of gap widths for linear arithmetic progressions and more
general nonlinear sequences
\cite{threegap,gapproblems,filtri,distinct,sqrts,logs}.  While these
studies cover the statistical distribution of the gap width
regardless of their position on the unit interval, the dependence of
the gap width on the position is a separate topic that is
presented in this article.

Consider the sequence of fractional parts of square roots of consecutive
integers,
\begin{equation}
  a_n=\{\sqrt{n}\},\quad n\in \mathbb{N},
\end{equation}
where $\{\}$ operator refers to the fractional part.
Taking more and more terms of the sequence, the unit interval is filled
asymptotically uniformly, but all the terms of the
sequence are either irrational or $0$. After $N$ terms, the average gap width is
obviously $1/N$, and the filling is uniform
enough that the maximum gap also converges to zero (no holes are
left). However, special behaviour is expected around rational points on the
unit interval; the gaps around them scale with an exponent slower than $N^{-1}$.

In this article, we go through derivation of asymptotic behaviour of rational
gaps for the sequence of square roots and its generalizations to higher-order
roots and subsequences. We show that the gap widths scale with a function of
rational argument that can be explicitly derived and is related to the Thomae's function.
We also investigate the geometric interpretation of the result for the sequence
of square roots.

\section{The gap function}

\subsection{Notation}

The modulo operator has the conventional meaning of describing
equivalence classes. We will use it to represent periodically
extended sets, which will enable a more condensed notation.
\begin{equation}
 \{a+bn \mid n \in\mathbb{N}\}=\{a \mod b\}.
\end{equation}

If two periods are in effect, the set can be reparameterized to cover a single period:
\begin{equation}
  \{a + b n \mod c \mid n \in \mathbb{N}\}=\{a \mod \gcd(b,c)\}.
  \label{eq:absorb}
\end{equation}
Element-wise multiplication of the set with a scalar multiplies both
the expression and the period (modulus),
\begin{equation}
  c\{ a\mod b \}=\{ac\mod bc \},
\end{equation}
so every common factor can be carried to the front.

Define the gap operator on the set $A$, centered at $x$:
\begin{equation}
  \gap_x(A)=\inf \{ y \mid y\in A \wedge y> x \}-\sup \{ y \mid y\in A \wedge y<x \}.
  \label{eq:gapdef}
\end{equation}
The center $x$ is assumed zero if omitted. In this homogeneous case, the following useful property holds for element-wise multiplication of a set with a scalar,
\begin{equation}
  \gap(tA)=|t|\gap(A).
  \label{eq:gapfactor}
\end{equation}

\subsection{Definition}

Define a set of fractional parts of square roots of first $N$ positive integers:

\begin{equation}
  f(N) = \{ \sqrt{n} \mod 1 \mid  n \in \mathbb{N} \wedge n\leq N \wedge \sqrt{n} \notin \mathbb{N} \}
\end{equation}

This is a finite set over the open unit interval
$(0,1)$ with cardinality $\#f(N)=N-\lfloor \sqrt{N}\rfloor$, 
periodically extended to include all numbers with the same fractional parts.

The main focus of this work is the following function that maps a real number to the width of
the gap in $f(N)$ around it:
\begin{equation}
  g(x,N) = \gap_x(f(N)).
\end{equation}

Because there are no rational numbers in the set $f(N)$, the gap widths around rational numbers have an interesting behaviour. Let us define the gap function as the limit
\begin{equation}
  G(x)=\lim_{N\to\infty}2\sqrt{N}g(x,N).
  \label{eq:Glim}
\end{equation}

For irrational values of $x$, the limit is expected to converge to $0$, as the average gap width falls as $N^{-1}$, which is faster than $N^{-1/2}$. The approximants to $G(x)$ at irrational $x$ thus fall off as $2/\sqrt{N}$. These ``background'' gaps follow an unusual distribution. Rather than having an exponential tail, as one would expect for a random distribution on a unit interval, the distribution has a $t^{-3}$ power-law asymptotic behaviour \cite{sqrts}.

In the following sections, we derive the closed form expression for the function $G(x)$ and its generalizations.

\subsection{Derivation}

Consider first the gaps on the open interval $(0,1)$, ignoring the singular case of the gap around the integer congruence class (the gap $G(0)$).

Between consecutive perfect squares, $\{\sqrt{n}\}$ runs across the unit interval. Up to the upper limit $N$, there are $\sqrt{N}$ passes. We will parameterize the running integer $n$ with two integer parameters,
\begin{align}
  n=k^2+m;&\quad m\in [1,2k] \\
  \label{eq:mdef}
  \sqrt{n}=k+x;&\quad x\in (0,1)
\end{align}
where $k$ is the integer part of the square root, $x$ is the fractional part, and $m$ is the ``remainder''. Together, this yields a quadratic equation for $x$,
\begin{equation}
  x^2+2kx-m=0.
  \label{eq:quadratic}
\end{equation}

Let $y=\frac{p}{q}$ be a target \emph{rational} number at which we want to evaluate the gap function ($p$ and $q>1$ are coprime). We must find the closest two irrational numbers that satisfy equation Eq.~\ref{eq:quadratic} and bracket our value $y$. For large enough $N$, the interval around $y$ is densely populated and the distance to the closest square root fractional part $x$ is small. By using the $y$ as an initial value, one iteration Newton's method on Eq.~\ref{eq:quadratic} approximates this distance, and becomes exact in the limit $N\to\infty$.
\begin{equation}
  \Delta(p/q,k,m)=x-y\approx\frac{-y^2-2ky+m}{2y+2k}=\frac{1}{2(k+y)}\left(-\frac{p^2+2kpq}{q^2}+m\right).
  \label{eq:newtstep}
\end{equation}
Among all the values this expression can assume for $n=k^2+m\leq N$, the smallest positive and largest negative value of $\Delta$, happen for consecutive $m$ and bracket a gap around zero, which approximates the gap $g(p/q,N)$. For each $k$, a large enough $m$ exists that makes the expression in the parentheses positive. Because $(p^2+2kpq)/q^2<2k+1$, such $m$ satisties the condition $m<2k+1$. Therefore we can safely drop this condition and allow $m$ to be any integer without affecting the gap width.
\begin{equation}
  g(p/q,N)\approx \gap(\Delta(p/q,k,m)\mid m,k\in \mathbb{N})
\end{equation}
The parenthesized part of Eq.~\ref{eq:newtstep} yields the same value for an infinite number of pairs $(k,m)$. The Diophantine equation $-p^2-2kpq+mq^2 = \const $ has periodic solutions: a pair $(k,m)$ gives the same value as $(k+tq/\gcd(2p,q) ,m+2tp/\gcd(2p,q) )$ for any integer multiplier $t$.

\begin{figure}
  \centering
  \includegraphics[width=\textwidth]{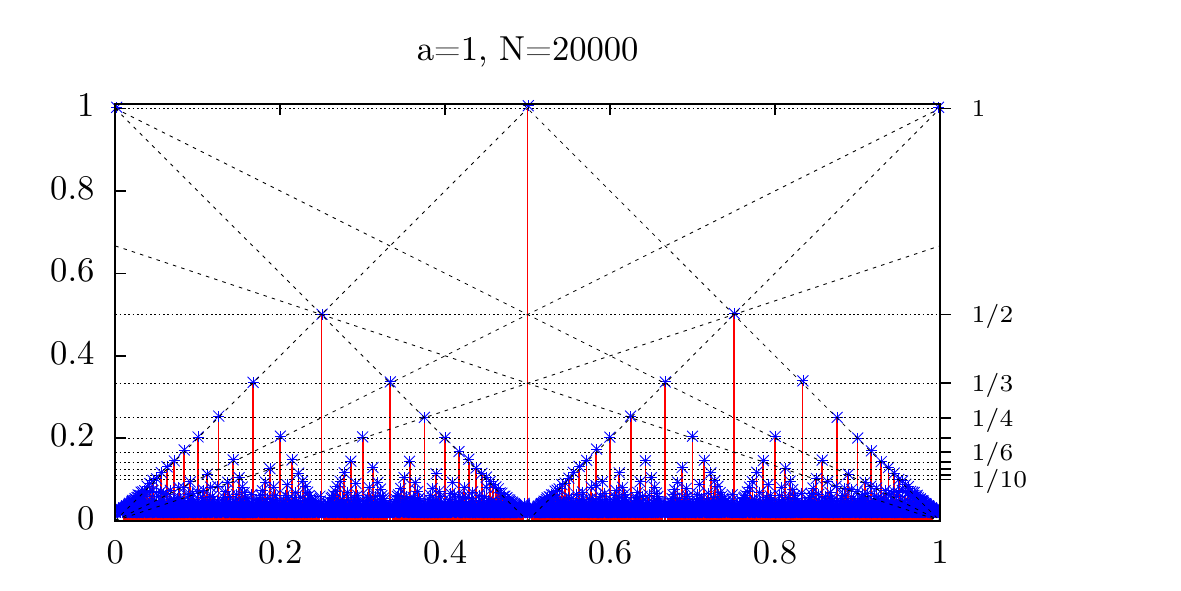}
  \caption{The function $G(x)$ and one of its numerical approximants computed at $N=20000$. The leftmost and rightmost data point represents the right and left half-gap, $G(0^+)$ and $G(1^-)$.}
  \label{fig:g2}
\end{figure}

\begin{figure}
  \centering
  \includegraphics[width=\textwidth]{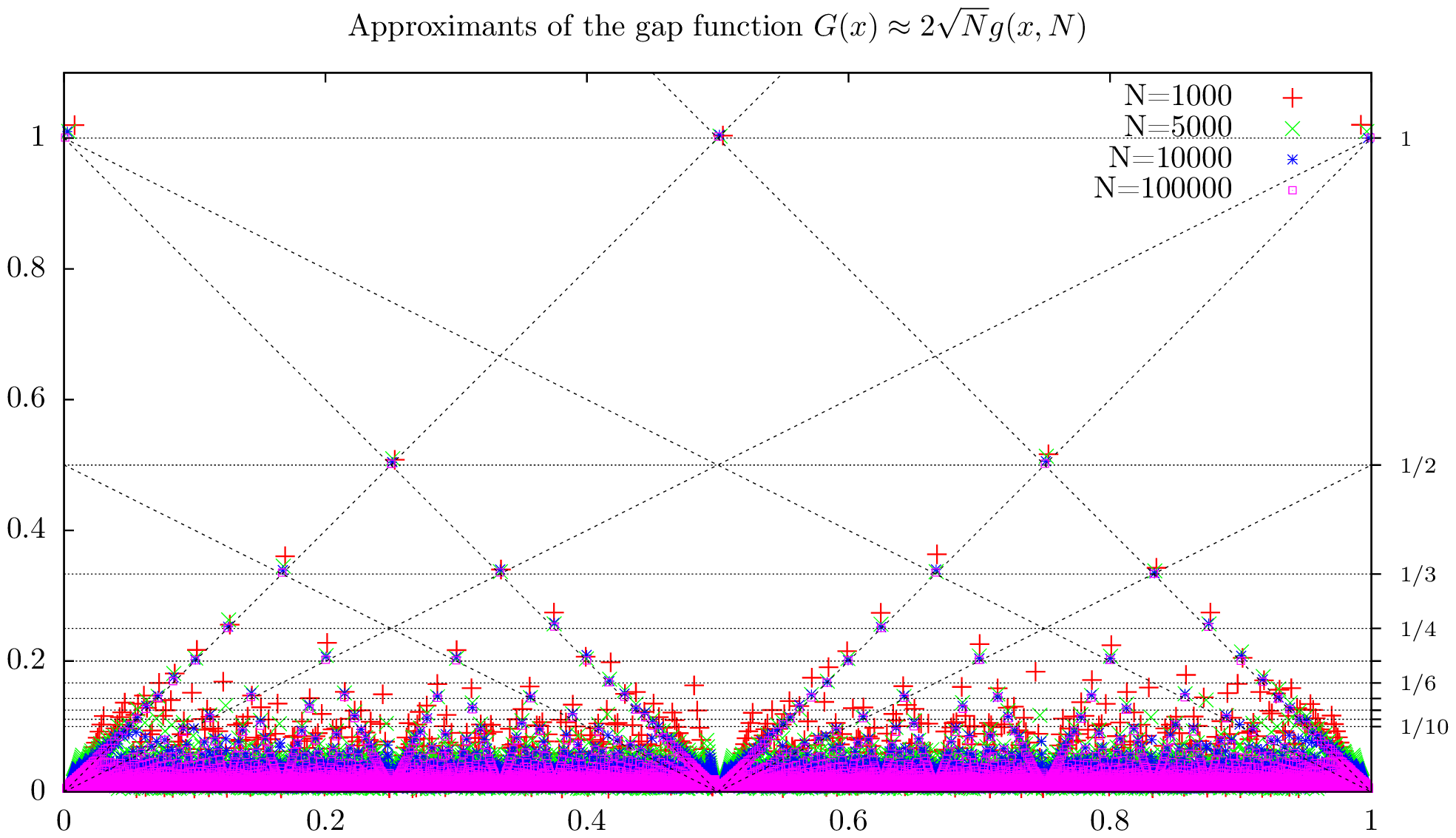}
  \caption{A series of numerical approximants to the gap function $G(x)$. Horizontal guidelines mark rational levels at $\frac{1}{q}$ where the limit values are supposed to lie, and a few diagonal lines with slopes $1$ and $2$ are used to highlight the collinearity of the values. Observe how the points converge to their limiting values, while the average background level is $2/\sqrt{N}$.}
  \label{fig:fig_gaps}
\end{figure}

For any given value of the parenthesized expression in Eq.~\ref{eq:newtstep}, the pairs $(k,m)$ that minimize the gap are those which minimize the pre-factor $\frac{1}{2(k+y)}$. The lowest values are reached with the highest $t$ allowed by the current upper bound $k^2+m\leq N$. This means that the pairs $(k,m)$ that actually gap the interval around $0$ have $k\in (\sqrt N-q/\gcd(2p,q),\sqrt N)$. In the limit $N\to \infty$, we have $q\ll \sqrt N$ and the prefactors tend to
\begin{equation}
  \frac{1}{2(k+y)}\asymp\frac{1}{2\sqrt N}.
\end{equation}
Note that the derivative in the denominator of the Newton's step only set the prefactor, which only depends on the $N$, and the
numerator (the parenthesised part of the expression Eq.~\ref{eq:newtstep}) by itself does not depend on $N$, so the limit in Eq.~\ref{eq:Glim} can be taken explicitly:
\begin{equation}
  G\left(\tfrac{p}{q}\right)= \gap\left\{\frac{p^2+2kpq}{q^2}-m\middle| m,k \in \mathbb{N} \right\}.
\end{equation}
where we utilized Eq.~\ref{eq:gapfactor} to carry the constant factor out and omit the negation. Furthermore, $m$ can be absorbed into the modulus by virtue of Eq.~\ref{eq:absorb}, and in the next step, so can $k$,
\begin{align}
  G\left(\tfrac{p}{q}\right)=& \gap\left\{\frac{p^2+2kpq}{q^2} \mod 1\middle| k \in \mathbb{N} \right\}=\\
  =& \frac{1}{q^2}\gap\left\{p^2+2kpq \mod q^2\middle| k \in \mathbb{N} \right\}=\\
  =& \frac{1}{q^2}\gap\left\{p^2\mod \gcd(2pq,q^2) \right\}=\\
  =& \frac{1}{q^2}\gap\left\{p^2\mod q\gcd(2,q) \right\}.
\end{align}
We used the fact $\gcd(p,q)=1$ to simplify the modulus. As the set we are observing is simply a periodic equally spaced lattice, all the gaps are equal to the period (modulus), and the offset $p^2$ plays no role in the result. Finally, we can write down the solution,
\begin{equation}
  G\left(\tfrac{p}{q}\right) = \frac{\gcd(2,q)}{q}=\frac{d}{q}
\end{equation}
where we marked $d=\gcd(2,q)\in \{1,2\}$. This expression can be written more explicitly as
\begin{equation}
  G(x)=
  \begin{cases}
    2 & x=0\\
    \frac{2}{q} &  x = \frac{p}{q}, q = 0 \mod 2\\
    \frac{1}{q} &  x = \frac{p}{q}, q = 1 \mod 2\\
    0 & x\notin \mathbb{Q}
  \end{cases},
  \label{eq:gx}
\end{equation}
and is depicted in Fig.~\ref{fig:g2} along with one of its approximants for a finite set generated by radicals up to $N=20000$.
A converging sequence of numerical approximants of this function is shown in Figure \ref{fig:fig_gaps}.

The gap function $G(x)$ is actually a familiar sight; it is a rescaled version of the Thomae's function $T(x)=G(2x)$ on the interval $(0,1)$ (see \cite{thomae,nechaev}, for instance).

\subsection{The gap around integers}
\label{sec:whole}

The gap $G(0)$ is is composed of the right half-gap $G(0^+)$ and left half-gap $G(0^-)=G(1^-)$ and as we excluded square roots of perfect squares from our set, the full gap is a sum of both half-gaps. The $G(0^+)$ gap is bounded from above by square roots of numbers right above perfect squares: $m=1$ and $p=0$. Equation \ref{eq:newtstep} reduces to $\frac{1}{2(k+y)}$, leading to $G(0^+)=1$. The left half-gap is obtained when $m=2k$ and $p=q$ ($y=1$), and is equal to $G(0^-)=1$, and together, $G(0)=2$.

\clearpage

\section{The generalized gap function}
\label{sec:general}

\subsection{Definition}

While the original gap function considered the square roots of all integers, we can generalize it to only include multiples of a chosen integer parameter $a$, falling back to the original definition when $a=1$.
\begin{equation}
  f^{(a)}(N) = \{ \sqrt{an}\mod 1 \mid   n \in \mathbb{N} \wedge an\leq N \wedge \sqrt{an} \notin \mathbb{N} \}
\end{equation}

\begin{equation}
  g^{(a)}(x,N) = \gap_x(f^{(a)}(N))
\end{equation}

\begin{equation}
  G^{(a)}(x)=\lim_{N\to\infty}2\sqrt{N}g^{(a)}(x,N)
\end{equation}

As the generalized set is only diluted (some terms are omitted), the gaps can
either stay equal, or become wider than before.

\subsection{Derivation}

In this section, we derive the function $G^{(a)}(x)$ using a similar
procedure used on the regular square root gaps.

A generalization to $\sqrt{an}$ does not affect the procedure up to
the Newton step performed in Eq.~\ref{eq:newtstep}, but it does
restrict the selection of $m$ for each $k$ differently, making them
coupled. To proceed with the reduction, we must find two integer
parameters that can be varied independently. The generalization of
Eq.~\ref{eq:mdef} is $m=an-k^2$, and substituting this relation, we
get a parameterization with $(k,n)$ which is indeed a pair of
independent positive integers. The generalized gap function is now
\begin{align}
  G^{(a)}\left(\tfrac{p}{q}\right)=& \gap\left\{\frac{p^2+2kpq+k^2q^2}{q^2}-an\middle | k,n\in\mathbb{N}\right\}=\\
  =& \frac{1}{q^2}\gap\left\{p^2+2kpq+k^2q^2 \mod aq^2\middle | k\in\mathbb{N}\right\}.
\end{align}

This may look similar to the previous case, but the quadratic term
$k^2$ makes the values in the set non-equally spaced and the offset
$p^2$ cannot be neglected.

The presence of $a$ and $q$ factors in the modulus suggests
decomposition $k=Ba+C$ with $C\in \mathbb{Z}_a$ and $B\in
\mathbb{Z}_q$.
\begin{equation}
  G^{(a)}\left(\tfrac{p}{q}\right)= \frac{1}{q^2}\gap\left\{p^2+2Bapq+2Cpq+q^2C^2 \mod aq^2 \middle | B\in \mathbb{Z}_q,C\in\mathbb{Z}_a\right\}.
\end{equation}
The parameter $B$ only appears linearly and rescales the modulus
according to Eq.~\ref{eq:absorb}. The new modulus is
$\gcd(2apq,aq^2)=aqd$.
\begin{align}
  G^{(a)}\left(\tfrac{p}{q}\right)=& \frac{1}{q^2}\gap\left\{p^2+2Cpq+q^2C^2 \mod aqd \middle | C\in\mathbb{Z}_a\right\}\\
  =& \frac{1}{q^2}\gap\left\{p^2+qd((2/d)Cp +(q/d)C^2) \mod aqd \middle | C\in\mathbb{Z}_a\right\}.
\end{align}
The first term can be forcibly split,
\begin{equation}
p^2=qd\lfloor p^2/qd \rfloor+\underbrace{p^2\mod qd}_{\delta},
\end{equation}
where $\mod$ here refers to the remainder after division, not the entire congruence class.

\begin{equation}
  G^{(a)}\left(\tfrac{p}{q}\right)= \frac{1}{q^2}\gap\left\{\delta+qd( \lfloor p^2/qd \rfloor +(2/d)Cp +(q/d)C^2) \mod aqd\middle | C\in\mathbb{Z}_a\right\}
\end{equation}
Because $\delta$ is smaller than the stride $qd$ of the second term, it is located in the same gap as zero, and can be omitted. Then, the common factor of the stride and the modulus can be carried outside.
\begin{equation}
  G^{(a)}\left(\tfrac{p}{q}\right)= \frac{d}{q}\gap\left\{\epsilon+\lfloor p^2/qd\rfloor+(2/d)Cp +(q/d)C^2 \mod a\middle | C\in\mathbb{Z}_a\right\}
  \label{eq:finalform}
\end{equation}
Note that in this form, the set could contain zero, which would make the gap definition ambiguous.
To amend this, we kept an infinitesimal displacement $\epsilon$ to specify that the zero value counts as an upper boundary of the gap. The expression is straight-forward to compute numerically, as it only involves finding the largest and the smallest value of a quadratic polynomial over $\mathbb{Z}_a$. The prefactor equals the expression for $a=1$, and the new $\gap$ factor tells us how much wider each gap is because of the skipped terms. In a special case when $q/d = 0 \mod a$, we are back to the linear case, and the gap stays the same width.

\begin{figure}
  \centering
  \includegraphics[width=0.75\textwidth]{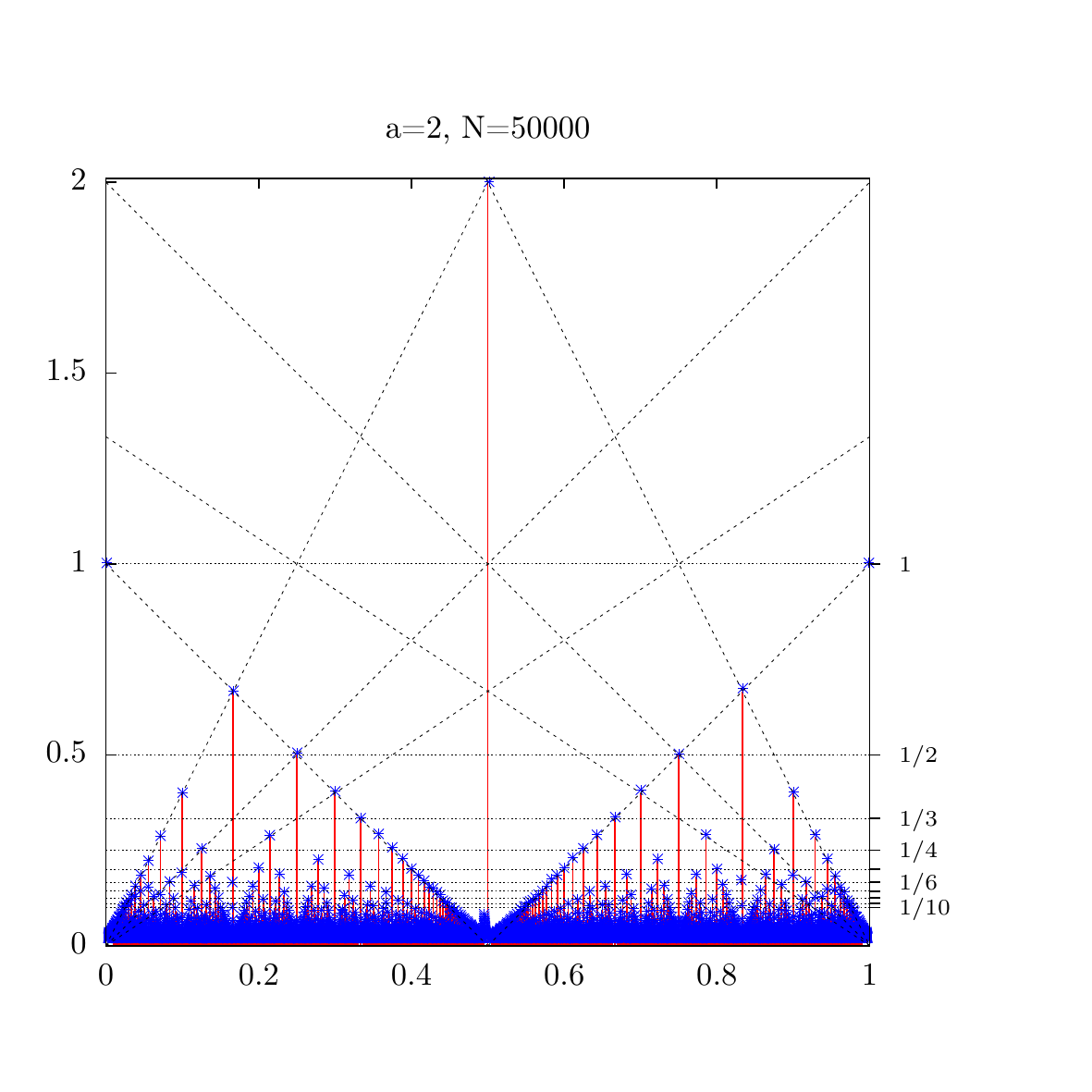}
  \caption{The function $G^{(2)}(x)$ and one of its approximants. Because for $a=2$, the quadratic residue term is suppressed, the pattern is still regular and recognizable, and divides the rational numbers into 3 cases based on the denominator modulo $4$ (Eq.~\ref{eq:case2}). Guidelines with slope $4$ are added to follow the wider gaps caused by the missing values.}
  \label{fig:g2a2}
\end{figure}

\begin{figure}
  \centering
  \includegraphics[width=0.75\textwidth]{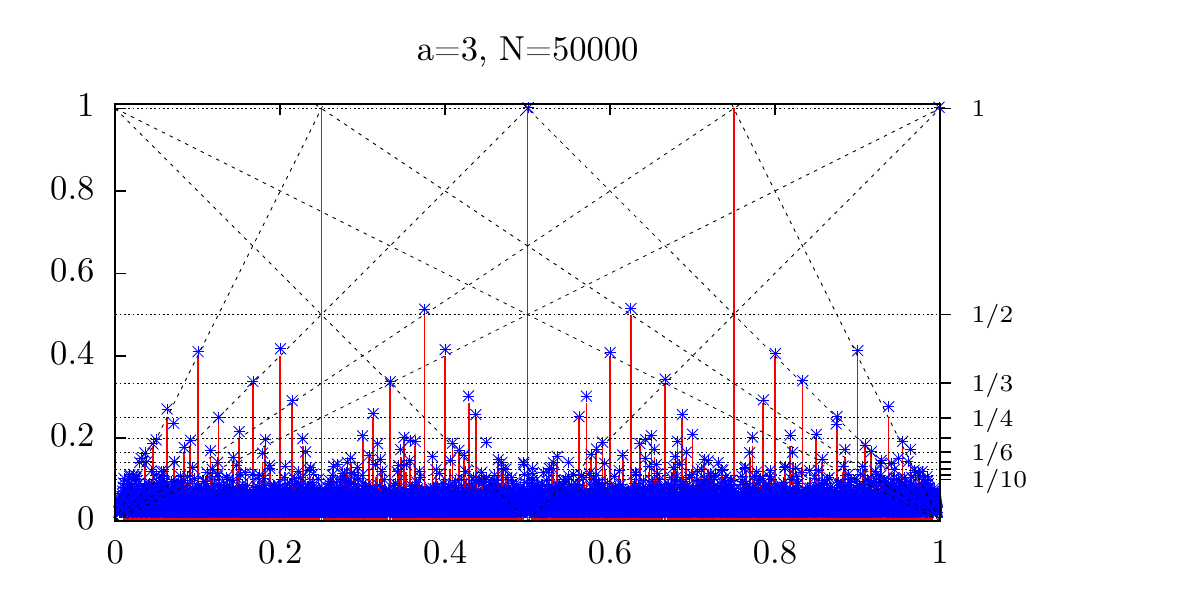}
  \caption{The function $G^{(3)}(x)$ and one of its approximants. This case has a more complex structure governed by the gaps in quadratic residues.}
  \label{fig:g2a3}
\end{figure}

\begin{figure}
  \centering
  \includegraphics[width=0.75\textwidth]{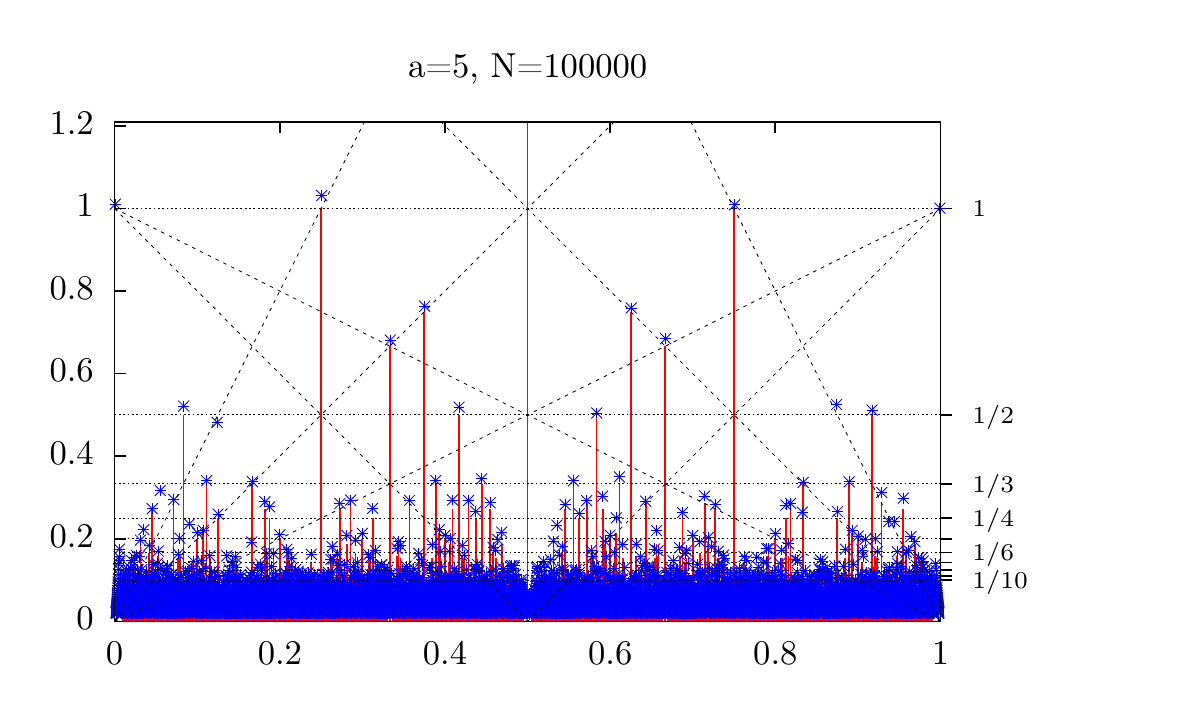}
  \caption{A higher-order example, $G^{(5)}(x)$. Observe how the spike heights are following a seemingly less organized sequence. The value $G^{(5)}(1/2)=6$ is clipped.}
  \label{fig:g2a5}
\end{figure}

Examples for the first few values of $a$ are shown in Figs.~\ref{fig:g2a2},\ \ref{fig:g2a3},\ \ref{fig:g2a5}.

The form $\sqrt{an}$ can be further generalized to $\sqrt{an+b}$ without much difficulty. The only difference is to substitute $p^2$ with $p^2-bq^2$ in Eq.~\ref{eq:finalform}. This leads to different gap factors for each $b$.

\subsection{The gap around integers}

The gap $G^{(a)}(0)$ is also affected by the condition $an=k^2+m$ which can be restated as $m=-k^2\mod a$. We can write it in the gap form,
\begin{equation}
  G^{(a)}(0)=\gap\left\{-k^2 \mod a\middle| k\in \mathbb{Z}_a \wedge k\neq 0 \right\}.
\end{equation}
The value $m=-1\mod a$ is a negative quadratic residue ($k=1\mod a$), but the $m=1\mod a$ may not be. We can simplify it to
\begin{equation}
  G^{(a)}(0)=\left(1+\inf\left\{-k^2 \mod a\middle| k\in \mathbb{Z}_a \wedge k\neq 0 \right\}\right).
\end{equation}

\subsection{The example for $a=2$}

If we set $a=2$, the gap function can still be reduced to an explicit form, because $C^2=C \mod 2$ and the quadratic term disappears. Equation \ref{eq:finalform} reduces to
\begin{equation}
    G^{(2)}\left(\tfrac{p}{q}\right)= \frac{d}{q}\gcd((2/d)p+(q/d),2)
\end{equation}
or, case-by-case,
\begin{equation}
  G^{(2)}\left(\tfrac{p}{q}\right)=\begin{cases}
  \frac{4}{q} & q=2 \mod 4 \\
  \frac{2}{q} & q=0 \mod 4 \\
  \frac{1}{q} & q=1 \mod 2.
  \end{cases}
  \label{eq:case2}
\end{equation}
This function and its approximants are shown in Fig.~\ref{fig:g2a2}.

\subsection{Geometric interpretation of the gap function}

The authors of Ref.~\cite{sqrts} study the probability distribution of
gap widths by observing a square lattice and how random lattice
translates overlap a given triangle. We take the same idea,
but instead of studying probability across all configurations, we 
retain the positional information and study projections of the lattice
onto a line.

\begin{figure}
  \centering
  \includegraphics[width=0.45\textwidth]{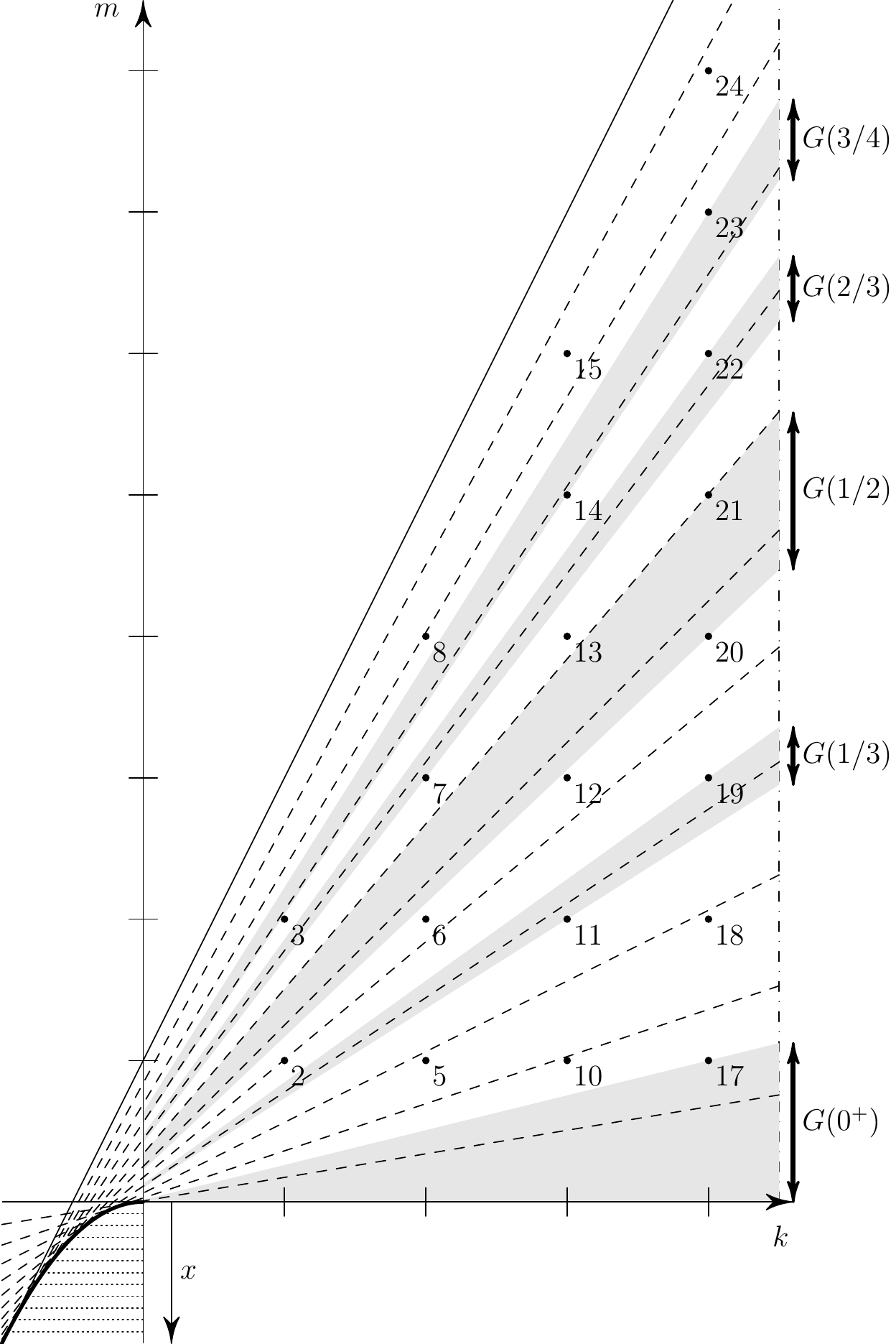}
  \includegraphics[width=0.45\textwidth]{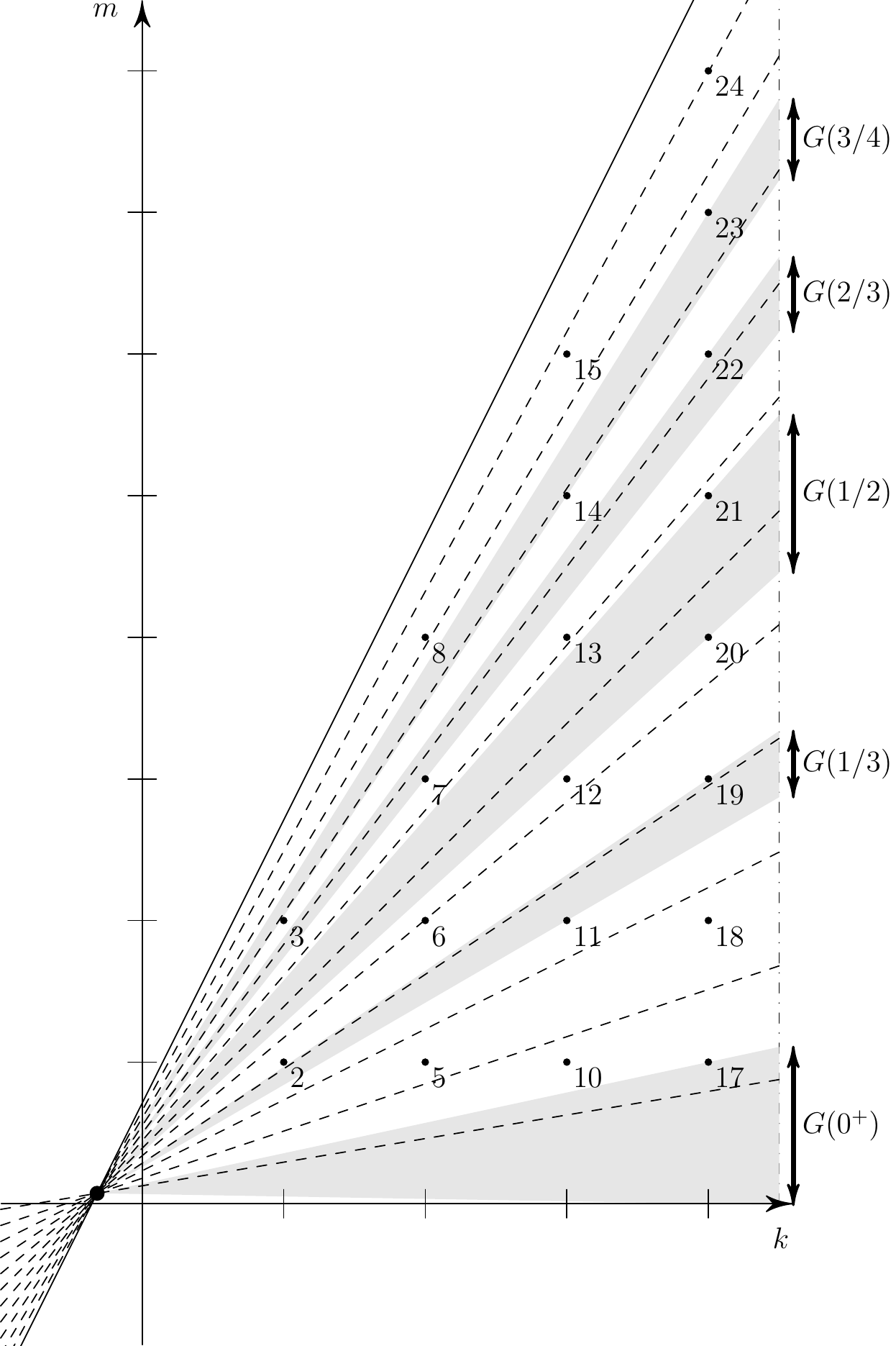}
  \label{fig:lattice}
  \caption{The Euclid's orchard representing the $(k,m)$ lattice, and the rays that project the lattice to the screen at the right. The ray pencils
    around rational slopes $2x$ pass unobstructed between the lattice points, their width equal in the limit to
    the gap function. {\sl Left:} The family of rays that corresponds to the fractional parts of square roots
    form a parabolic caustic with the vertical component of the point of tangency equal to $-x$. {\sl Right:}
    Asymptotically, any point with irrational coordinates yields the same light pattern on the screen. Thus, the gap function is the illumination pattern of a single point source projected through the Euclid's orchard.
  }
\end{figure}

The defining equation for the elements $x$ of the sequence of fractional
parts of square roots (Eq.~\ref{eq:quadratic}) can be seen as a function
in the $(k,m)$ plane,
\begin{equation}
  m(k)=x^2+2kx,
  \label{eq:mgeom}
\end{equation}
which describes a family of rays tangent to the caustic $m=-k^2$
parameterized by $x$, with the point of tangency at $(-x,-x^2)$
(Figure 2 in \cite{sqrts} and \ref{fig:lattice}).

Consider the Euclid's orchard: a lattice of integer points $(k,m)$ in two dimensions.
The gap $g(x,N)$ measures the width of the interval where no $x$
solves Eq.~\ref{eq:mgeom} for integer pairs $(k,m)$ in the range
$k<\sqrt{N}$. Set a vertical screen (a line) at $k_{\text{max}}=\sqrt{N}$.
A distance swept by the ray on the screen without casting a shadow of a lattice point,
simplifies to
\begin{equation}
  x_2^2+2k_{\text{max}}x_2-x_1^2-2k_{\text{max}}x_1=2(x_2-x_1)(k+x)=2(k+x)g(x,N)\to G(x),
\end{equation}
where the mean value $x=(x_1+x_2)/2$ was used to represent each
interval. Notice that the expression is exactly the parenthesized part
of Eq.~\ref{eq:newtstep} and in the limit becomes the gap function.

With our specific choice of the ray family, the gap function directly
measures the length of contiguously illuminated parts of the
screen. Recall, however, that the $x^2$ term in Eq.~\ref{eq:mgeom}
corresponds to the $p^2$ term in Eq.~\ref{eq:newtstep}, the value of
which does not affect the gap. The only exception is the singular case
$p=0 \mod q$ when a term of the sequence hits a rational point instead
of leaving a gap around it. In our case, the perfect squares which hit
the value $x=0$ were omitted from the sequence and the gaps around
the integers were treated separately in section \ref{sec:whole}.

Instead of $x^2$, the $m$-intercept $c(x)$ of the rays can be almost
any function of $x$,
\begin{equation}
  m(k)=c(x)+2kx,
\end{equation}
and still produce the same illumination pattern in the limit.  The
only condition is to avoid the singular case $c(p/q) = 0 \mod 1$.
Such cases again correspond to rational values in the
sequence. Geometrically, these are the rays that pass through an
infinite set of lattice points, casting a single point shadow instead
of leaving an illuminated gap around the point.

Perhaps the simplest candidates for the $m$-intercept are linear,
$c(x)=c_1+c_2 x$, forming a pencil of rays passing through a single
focus instead of a caustic. If $c_1$ and $c_2$ are in a rational
proportion, we get some rational terms in the sequence. In the most
singular case $c(x)=0$, the entire sequence is rational and instead of
gaps, the point shadows form an image of the rational set on the
screen. Recall that perspective projection of an orchard of unit height trees
produces a Thomae's function, but in perceived heights, not in the lengths of
illuminated gaps.

For all other choices of $c_{1,2}$, the projection of the lattice onto
the screen leaves gaps that converge toward the gap function we
defined in the beginning, though the continued fraction expansion
determines how quickly the limit converges for the closest rational
approximations of the ratio $c_1/c_2$. The invariance to the
translation of rays is related to the ergodicity that enabled
calculation of the gap distribution in Ref.~\cite{sqrts}.

In short, the gap function tells us the lengths of the portions of a
screen illuminated by a point source located at irrational
coordinates, shining through a square lattice -- an Euler's
orchard. The generalized gap function corresponds to the pattern
obtained by removing some of the lattice points so that only the
multiples of a certain integer cast shadows.

\clearpage

\section{The higher order root gaps}

\subsection{Definition and derivation}

Another straight-forward generalization is the sequence $\{\sqrt[\alpha]{n}\}$. We can proceed as before, by defining
\begin{equation}
  n=k^\alpha+m;\quad m\in [1,(k+1)^\alpha-k^\alpha-1]
  \label{eq:mdefalpha}
\end{equation}
and $\sqrt[\alpha]{n}=k+x$. Instead of a quadratic equation, we obtain a polynomial
\begin{equation}
  (x+k)^\alpha-k^\alpha-m=0.
  \label{eq:ndratic}
\end{equation}
We must check if the Newton step procedure is still valid for a higher order polynomial we are dealing with now. It turns out that for $x\in (0,1)$, the method converges, as the above polynomial is monotonously increasing from negative to positive on this domain with a single primitive positive root. The delta step yields
\begin{equation}
  \Delta(p/q,k,m)=\frac{1}{\alpha(k^{\alpha-1}+\mathcal{O}(k^{\alpha-2}))}
  \left(-\frac{\sum_{i=0}^{\alpha-1} {\alpha \choose i} p^{\alpha-i}(kq)^{i}}{q^{\alpha}}+m\right).
\end{equation}

The prefactor shows the power law scaling that has to be adjusted for the $\alpha$-th order gap function to exist in the limit:
\begin{align}
  G_\alpha (p/q)=&\lim_{N\to\infty}\alpha N^{(\alpha-1)/\alpha} \gap(\Delta(p/q,k,m) \mid m,k\in\mathbb{N})=\\
  =&\frac{1}{q^\alpha}\gap\left\{\sum_{i=0}^{\alpha-1} {\alpha \choose i} p^{\alpha-i}(kq)^{i} \mod q^\alpha \middle | k\in\mathbb{N} \right\}.
\end{align}
the prefactor scaling is getting closer and closer to the average gap distribution, which goes as $1/N$. This means that it takes a larger $N$ to suppress the background compared to the gap function, and the limit converges very slowly for $\alpha>2$.

We learned from the $G^{(a)}$ generalization that it is convenient to express $k$ in base $q$:
\begin{equation}
  k=\sum_{j=0}^{\alpha-2} C_j q^{j}
\end{equation}
where $C_j\in \mathbb{Z}_q$. The term $C_{\alpha-1} q^{\alpha-1}$ is not needed because $q C_{\alpha-1}q^{\alpha-1} = 0\mod q^{\alpha}$. Let us also split off the first two terms of the sum over $i$:
\begin{equation}
  G_\alpha (p/q)=\frac{1}{q^\alpha}\gap\left\{p^\alpha+\alpha p^{\alpha-1}q\sum_{j=0}^{\alpha-2}C_jq^j+ \sum_{i=2}^{\alpha-1} {\alpha \choose i} p^{\alpha-i}q^{i}\left[\sum_{j=0}^{\alpha-2} C_j q^{j}\right]^{i} \mod q^\alpha \middle | C_j\in\mathbb{Z}_q \right\}.
  \label{eq:Galpha_unreduced}
\end{equation}
In the last term, the prefactor $q^{i\geq 2}$ ensures that $C_{\alpha-2}$ contribution to the sum is a multiple of $q^\alpha$. The second term, which is linear in $C_j$, is the only surviving term that contains $C_{\alpha-2}$. The parameter $C_{\alpha-2}$ can be absorbed, rescaling the modulus to
\begin{equation}
  \gcd(\alpha p^{\alpha-1} q^{\alpha -1},q^\alpha)=\gcd(\alpha,q)q^{\alpha-1}=d_1 q^{\alpha-1}.
\end{equation}
The reduced expression becomes
\begin{equation}
  G_\alpha (p/q)=\frac{1}{q^\alpha}\gap\left\{p^\alpha+\alpha p^{\alpha-1}q\sum_{j=0}^{\alpha-3}C_jq^j+ \sum_{i=2}^{\alpha-2} {\alpha \choose i} p^{\alpha-i}q^{i}\left[\sum_{j=0}^{\alpha-3} C_j q^{j}\right]^{i} \mod d_1q^{\alpha-1} \middle | C_j\in\mathbb{Z}_q \right\}.
\end{equation}
This step removed the highermost $C_j$ coefficient and reduced the powers of $q$ by one. Now the next biggest term is the one that only stands in the linear term, so the coefficient absorption can be performed recursively until the nonlinear part vanishes completely.

This procedure eventually yields a linear expression in the last remaining coefficient $C_1$. Because each next modulus divides the previous, all the steps of the recursion can be condensed into looking up the modulus induced by the term $\alpha p^{\alpha-1} C_1 q$ in the unreduced expression Eq.~\ref{eq:Galpha_unreduced}:
\begin{equation}
  \gcd(\alpha p^{\alpha-1} q, q^\alpha)=q\gcd(\alpha,q^{\alpha-1})=qd.
\end{equation}
The constant term $p^\alpha$ can again be skipped and the higher order gap function written out as
\begin{equation}
  G_\alpha(p/q)=\frac{\gcd(\alpha,q^{\alpha-1})}{q^{\alpha-1}}=\frac{d}{q^{\alpha-1}}.
\end{equation}

If $\alpha$ is prime, then $d=\gcd(\alpha,q)$ and can be either $1$ or $\alpha$.
\subsection{Additional generalizations}

The higher-order gap functions can also be generalized to use a
diluted set of integers, $\sqrt[\alpha]{an}$. The procedure is a
similar to the square root case. The solution involves finding a gap
in a polynomial expression of order $\alpha$ modulo $a$.

Another possible generalization is to introduce rational exponents,
i.e. $\sqrt[\alpha]{n^\beta}$. In this case, the resulting modular
expression is nonlinear in both independent parameters, which makes
reduction less trivial.

\subsection{Convergence rate}

The analytical result for the gap function is really difficult to demonstrate numerically even for $\alpha=3$. First of all, the $1/N$ average gap size means that the background noise in the numerical approximation scales as $3N^{(3-1)/3}(1/N)=3N^{-1/3}$. Assuming these background gaps are distributed exponentially \cite{sqrts}, we estimate the minimal $N$ for which there are only a couple of outliers above a level $\epsilon$:
\begin{equation}
  Ne^{-3\epsilon N^{1/3}}\sim 1.
\end{equation}
As the gap function for $\alpha=3$ has smaller values compared to the square root gap function, $\epsilon\sim 1/3$ is barely enough to keep the tallest spikes ($q=3$) level with the farthest outliers. This means that it's impossible to see any kind of signal below $N\sim 2\cdot 10^6$, and the first spikes can begin to grow out of the noise around $N\sim 4\cdot 10^8$ to reach the third tallest spike at $q=6$. This is drastically worse than the $\alpha=2$ case where $N\sim 1000$ was already enough to plot a decent graph (Fig.~\ref{fig:fig_gaps}).

For $\alpha = 4$, values of order $N\sim 10^{14}$ would be needed to resolve the first hints of the tallest spikes, and $N\sim 10^{18}$ to get anything useful, which makes naive numerical demonstration impractical.

\section{Conclusion}

With the exception of integer solutions, expressions involving roots
of integers cannot assume rational values. An immediate consequence is
that rational points on the unit interval have a special role as a
repulsor for the fractional parts of the sequence. In this report, we
worked out how wide the gaps are around each rational number in the
form of an exact analytical expression. The relative gap widths, which
are, incidentally, also in rational proportions, tell us how well can
a fraction can be approximated by an irrational number (a root), not
unlike how the continued fraction determines how well an
\emph{irrational number} can be approximated with a fraction. In a
sense, fractions with a lower denominator tend to be ``less
irrational'', creating a wider gap in the sequence. The methodology
used here can be extended to a wide variety of sequences, albeit the
last step may not have a closed form expression.

\end{document}